\documentclass[a4paper,12pt]{amsart}
\usepackage{amsmath,amssymb,amsthm,ascmac,fancybox}   
\usepackage{bm}                
\usepackage{tabularx}          
\usepackage[dvipdfm]{graphicx} 
\usepackage{comment} 
\allowdisplaybreaks[3]
\newcommand{\sa}   [1]{\mathcal{#1}}                       
\newcommand{\set}  [2]{\left\{#1;#2\right\}}               
  
\newcommand{\narabe}  [1]{\begin{smallmatrix} #1 \end{smallmatrix}}   
\newcommand{\norm} [2]{\left\| #1 \right\|_{#2}}               
\newcommand{\ab}  [1]{\left| #1\right|} 				
\newcommand{\qw}   [0]{\par{\noindent}}                    
\newcommand{\ve}{\varepsilon}

\newcommand{\N}{\mathbb{N}}			
\newcommand{\R}{\mathbb{R}}			

\newcommand{\sign}{\mathrm{sign}\,}  	
\newcommand{\esup}{\mathrm{ess}\sup}	
\newcommand{\einf}{\mathrm{ess}\inf}	
\newcommand{\alm}{\mathrm{a.e.}\,}		
\newcommand{\kko}  [1]{\left( #1\right)}  

\newcommand{\ds}{\displaystyle}

\title{inequalities for derivative of polynomial approximation with Erd{\H{o}}s type weight}
\author{Kentaro\ Itoh, Ryozi\ Sakai and Noriaki\ Suzuki}

\begin{document}
\

\maketitle

\section{{Introduction}}

Let $\R=(-\infty,\infty)$. We consider an exponential weight
\[ w(x)=\exp(-Q(x)) \]
on $\R$, where $Q$ is an even and nonnegative function on $\R$. Throughout
this paper we always assume that w belongs to a relevant class $\sa{F}_\lambda(C^2 +)$ (see Section 2). We consider a function $T = T_{w}$ defined by
\[ T(x):=\frac{xQ'(x)}{Q(x)},~~x\neq 0 \leqno{(1.1)} \]
is very important. We call $w$ a Freud-type weight if $T$ is bounded, and otherwise, $w$ is
called an Erd{\"{o}}s-type weight. 

For $x >0$, the Mhaskar-Rakhmanov-Saff number (MRS number) $a_x = a_{x}(w)$ of $w=\exp(-Q)$ is defined by a positive root of the equation
\[ x=\frac{2}{\pi}\int_0^1\frac{a_xuQ^\prime(a_xu)}{(1-u^2)^{1/2}}du. \leqno{(1.2)} \]
When $w = \exp(-Q) \in \sa{F}(C^2+)$, $Q'$ is positive and increasing on $(0, \infty)$, so that
$$
\lim_{x \to \infty} a_x = \infty  \ \ \mbox{and} \ \ 
\lim_{x \to \infty} \frac{a_x}{x} = 0 \leqno(1.3)
$$
hold. Note that those convergences are all monotonically.

\vspace{2ex}


\vspace{1ex}



\vspace{1ex}

For $1\le p\le\infty$, we denote by $L^p(I)$ the usual $L^p$ space on an interval $I$ in  $\R$. For $fw\in L^p(\R)\ (1\le p\le\infty)$, the degree of weighted polynomial approximation $E_{p,n}(w,f)$ is defined by
\[ E_{p,n}(w;f):=\inf_{P\in\sa{P}_n}\|w(f-P)\|_{L^p(\R)}, \]
where $\sa{P}_n$ is the class of all polynomials with degree $n$ at most.

\vspace{1ex}

In [1], H. N. Mhaskar obtained the following result. Let $w=\exp(-Q)$ be a Freud-type
weight, and let $1\le p\le \infty$. Let $f$ be absolutely continuous with $f^\prime w\in L^p(\R)$. There exists a constant $C>0$ such that for $n\in\N$ and $P\in\sa{P}_n$, if
\[ \|w(f-P)\|_{L^p(\R)}\le \eta \]
holds for $\eta >0$, then
\[ \|w(f^\prime-P^\prime)\|_{L^p(\R)}\le CE_{p,n-1}(w;f^\prime)+\frac{n}{a_n}\eta. \]
In particular, if we take $\eta=2E_{p,n}(w,f)$,
then 
\[ \|w(f^\prime-P^\prime)\|_{L^p(\R)}\le CE_{p,n-1}(w;f^\prime) \]
holds (see [1, Theorem 4.1.7]).

\vspace{1ex}

Now, we will obtain a counterpart of the above theorem for an Erd{\H{o}}s-type weight. The following result is main theorem in this paper.

\vspace{1ex}

\qw{\bf Theorem 1.1.} Let $w=\exp(-Q)$ be a Erd{\H{o}}s-type
weight. Suppose that
\[ T(a_n)\le c\kko{\frac{n}{a_n}}^{2/3} \leqno(1.4) \]
holds some constant $c\ge 0$. For any $C_1\ge 1$, $1\le p\le \infty$ and absolutely continuous function $f$ with $f^\prime w\in L^p(\R)$, if
\[ \|(f-P)w\|_{L^p(\R)}\le C_1E_{p,n}(w;f) , \leqno(1.5)  \]
holds for some $n\in\N$ and $P\in\sa{P}_n$, then there exists constant $C_2\ge 1$
\[ \|w(f^\prime-P^\prime)\|_{L^p(\R)}\le C_2(1+C_1)T^{3/4}(a_n)E_{p,n-1}(w;f^\prime) \leqno(1.6) \]
holds.

\vspace{1ex}

In particular, when we denote the $n$-degree best approximation polynomial of $f$ with the weight $w$ by $P_{f,n}$, that is,
\[ E_{p,n}(w;f)=\|w(f-P)\|_{L^p}(\R), \leqno(1.7) \]
then there exists a constant $C>0$ such that
\[ \|w(f^\prime-P^\prime)\|_{L^p(\R)}\le CT^{3/4}(a_n)E_{p,n-1}(w;f^\prime). \leqno(1.8) \]
Here if $p=\infty$, then we suppose that $f^\prime$ is continuous and $w(x)f^\prime(x)\to 0$ as $|x|\to \infty$. When $f^\prime$ is continuous and $w(x)f^\prime(x)\to 0$ as $|x|\to \infty$  we see $w(x)f(x)\to 0$ as $|x|\to \infty$. In fact, it follows from Lemma 2. 11.
\vspace{1ex}
\mbox{}\\
Using (1.7) and (1.8), we have the following corollary.\\
\vspace{1ex}\\
\qw{\bf Corollary 1.2.} Let $w\in\sa{F}(C^{2}+)$ satisfy the condition (1.4). Let $1\le p\le\infty$, for integer $r\ge 1$ and $n\ge 1$, there exists a constant $C_1>0$, and $f^{r-1}$ be absolutely continuous in each compact interval  with $f^{(r)} w\in L^p(\R)$, where for $p=\infty$, we suppose that $f^{(r)}$ is continuous and $w(x)f^{(r)}(x)\to 0$ as $|x|\to \infty$ (then we see $w(x)f^{(j)}(x)\to 0$ as $|x|\to \infty$, for $j=0,1,\cdots,r$). If for $P\in\sa{P}_n$, we have
\[ \|w(f-P)\|_{L^p(\R)}\le C_1T^{1/4}(a_n)E_{p,n}(w;f), \]
then there exists a constant $C_2>0$ such that
\[ \|w(f^{(j)}-P^{(j)})\|_{L^p(\R)}\le C_2T^{(2j+1)/4}(a_n)E_{p,n-j}(w;f^{(j)}). \]

\section{{class $\sa{F}(C^2+)$ and some fundamental notations}}

We say that an exponential weight $w=\exp(-Q)$ belongs to class $\sa{F}(C^2+)$, when $Q:\R\to[0,\infty)$ is a continuous and even function and satisfies the following conditions:\\
$~~$(a)$~~Q^\prime(x)$ is continuous in $\R$ and  $Q(0)=0$.\\
$~~$(b)$~~Q^{\prime\prime}(x)$ exists and is positive in $\R\setminus\{0\}$.\\
$~~$(c)$\ds~~\lim_{x\to\infty}Q(x)=\infty.$\\
$~~$(d)$~~$The function 
\[ T(x):=\frac{xQ'(x)}{Q(x)},~~x\neq 0 \leqno{(1.1)} \]
 is quasi-increasing in $(0,\infty)$(i.e. there exists $C>1$ such that $T(x) \le C T(y)$ whenever $0<x<y$), and there exists $\Lambda\in\R$ such that
\begin{align*} 
T(x)\ge \Lambda>1,~~x\in \R\setminus\{0\}.
\end{align*}
$~~$(e)$~~$There exists $C>1$ such that
$$
\frac{Q''(x)}{|Q'(x)|}\le C\frac{|Q'(x)|}{Q(x)},~~\alm~x\in \R.
$$
Moreover, if there also exist a compact subinterval $J (\ni 0)$ of $\R$, and $C>1$ such that
\begin{align*} 
 C \frac{Q''(x)}{|Q'(x)|}\ge \frac{|Q'(x)|}{Q(x)},~~\alm~x\in \R\setminus J.
\end{align*}

\vspace{2ex}

Let $w\in\sa{F}(C^{2}+)$ and $\lambda >0$. We write $w \in \sa{F}_\lambda(C^2 +)$ if  there exist $K>1$ and $C>1$ such that for all $|x| \geq K$, 
$$ \frac{|Q^\prime(x)|}{Q(x)^\lambda}\le C 
$$
holds.  

A typical example of Freud-type weight is $w(x)=\exp(-|x|^\alpha)$ with $\alpha>1$. 
It belongs to $\sa{F}_\lambda(C^{2}+)$ for $\lambda =1$. 
For 
$u\ge 0,\ \alpha>0$ with $\alpha+u>1$ and $l\in\N$, we set 
\[ Q(x):=|x|^u(\exp_l(|x|^\alpha)-\exp_l(0)), \]
where 
\[ \exp_l(x) :=\exp(\exp(\exp(\cdots(\exp x)))) (l-\mbox{times}).\] 
Then $w=\exp(-Q(x))$ is an Erd{\"{o}}s-type weight, which belongs to $\sa{F}_\lambda(C^{2}+)$ for $\lambda >1$.

\vspace{2ex}

Throughout this paper we always assume that $w$ belongs to a class $\sa{F}(C^2+)$. Let $\{p_n\}$ be orthogonal polynomials for a weight $w$, that is, $p_n$ is the polynomial of degree $n$ such that
$$ \int_\R p_n(x)p_m(x)w^2(x)dx = \delta_{mn}. $$
For a function $f$ with $fw\in L^p(\R)$, we set
\[ s_n(f)(x):=\sum_{k=0}^{n-1}b_k(f)p_{k}(x) ~~\mbox{where}~~b_k(f)=\int_\R f(t)p_k(t)w^2(t)dt \]
for $n\in\N$ (the partial sum of Fourier series). Also it is given by
\[ s_n(f)(x)=\int_\R K_n(x,t)f(t)w^2(t)dt,  \]
where
\[ K_n(x,t) = \sum_{k=0}^{n-1}p_k(x)p_k(t). \]
It is known that by the Cristoffel-Darboux formula
\[ K_n(x,t) = \frac{\gamma_{n-1}}{\gamma_n}\frac{p_n(x)p_{n-1}(t)-p_n(t)p_{n-1}(x)}{x-t} \]
holds, where $\gamma_n$ and $\gamma_{n-1}$ are the leading coefficients of $p_n$ and $p_{n-1}$, respectively.
Then 
\[ a_n \sim \frac{\gamma_{n-1}}{\gamma_n}. \]
The de la Vall\'{e}e Poussin mean $v_n(f)$ of $f$ is defined by
\[ v_n(f)(x):=\frac{1}{n}\sum_{j=n+1}^{2n}s_j(f)(x). \]
The Christoffel function $\lambda_{n}(x) = \lambda_n(w,x)$ is defined by
\[ \lambda_n(x):=\frac{1}{K_n(x,x)}=\kko{\sum_{k=0}^{n-1}p_k(x)^{2}}^{-1}. \]
Then 
$$
\lambda_n (x) = \inf_{P\in\sa{P}_{n-1}}\frac{1}{P(x)^{2}} \int_\R |P(t)w(t)|^2dt.
$$
holds on $\R$. We denote the zeros of the orthonormal polynomial $p_n(x)$ by $x_{n,n}<x_{n-1,n}<\dots<x_{1,n}$. And then we define Christoffel numbers $\lambda_{k,n},\ k=1,2,\dots ,n$ such as $\lambda_{k,n}:=\lambda_n(x_{k,n})$. We also use the following function for $u>0$:
$$
\varphi_u(x) :=
\begin{cases}
       \displaystyle \frac{a_u}{u}\frac{1-\frac{|x|}{a_{2u}}}{\sqrt{1-\frac{|x|}{a_u}+\delta_u}},& |x|\le a_u, 
       \\[5ex]
         \varphi_u(a_u),& a_u<|x|, 
\end{cases}
$$
where $\delta_u :=\{uT(a_u)\}^{-2/3}$.

The facts described below are very important.
\vspace{1ex}\\
\qw{\bf Lemma 2.1.} (1) [2, Lamma 3.5 (3.27)-(3.29)] For fixed $L>0$ and $t>0$,
\[ a_{Lt}\sim a_t,\ \ \ \ \ \ \ \ T(a_{Lt})\sim T(a_t), \]
and for $j=0,1$
\[ Q^{(j)}(a_{Lt}) \sim Q^{(j)}(a_t). \]
(2) [2, Lamma 3.4 (3.17), (3.18)] For $t>0$,
\[ Q^\prime(a_t)\sim\frac{t\sqrt{T(a_t)}}{a_t},\ \ \ \ \ \ \ \ Q(a_t)\sim\frac{t}{\sqrt{T(a_t)}}. \]
(3) [2, Lamma 3.11] For fixed $L>0$ and $t>0$,
\[ \ab{1-\frac{a_{Lt}}{a_t}} \sim \frac{1}{T(a_t)}. \]
\vspace{1ex}\\
\qw{\bf Lemma 2.2.} (a) [2, Corollary 13.4 (a) (12.20)] Uniformly for $n\ge 1$, $1\le k\le n-1$,
\[ x_{k,n}-x_{k-1,n}\sim\varphi_n(x_{k,n}),\ \ \ \ \ \ \ \ \varphi_n(x_{k,n})\sim\varphi_n(x_{k+1,n}). \]
(b) [3, Lemma 3.4 (d)] Let $\max{|x_{k,n}|, |x_{k+1,n}|}\le a_{n/2}$,
\[ w(x_{k,n})\sim w(x_{k+1,n})\sim w(x),\ \ \ \ (x_{k+1,n}\le x\le x_{k,n}). \]
So, for given $C>0$ and $|x|\le a_{n/3}$, if $|x-x_{k,n}|\le C\varphi_n(x)$,
\[ w(x)\sim w(x_{k,n}). \]
\vspace{1ex}\\
\qw{\bf Lemma 2.3.} [2, Theorem 9.3 (c)] Let $1\le p\le \infty$ and $L>0$. Then uniformly for $n\ge 1$ and $|x|\le a_n(1+L\eta_n)$,
\[ \lambda_n(x)\sim\varphi_n(x)w^p(x). \]
Moreover, uniformly for $x\in\R$,
\[ \lambda_n(x)\ge C\varphi_n(x)w^p(x). \]
\vspace{1ex}\\

For any nonzero real valued functions $f(x)$ and $g(x)$, we write $f(x)\sim g(x)$ if there exist positive constants $C_1$ and $C_2$ independent of $x$, such that $C_1 g(x)\le f(x)\le C_2g(x)$ holds for all $x$. Similarly, for any positive numbers $\{c_n\}_{n=1}^\infty$ and $\{d_n\}_{n=1}^\infty$, we write $c_n\sim d_n$.

Throughout this paper, $C_1, C_2,\dots$ denote positive constants independent of $n, x, t$ or polynomials $P_n(x)$. The same symbol does not necessarily denote the same constant in different occurrences.\ We consider the numbers $p$ and $q$ are  $1\le p\le \infty$ with $1/p+1/q=1$, if $p=\infty$ then $q=1$ and  $p=1$ then $q=\infty$. \\

\section{{Lemmas}}

We prepare some important lemmas.

\qw{\bf Lemma 3.1.} [cf. 1, Appendeix A.1.1.] Let $fw\in L^p(\R)$, then
\[ E_{p,0}(w;f)=\sup_{\narabe{ \|gw\|_{L^q(\R)}\le 1 \\ \int_{\R}g(t)w^2(t)dt=0 }}\int_{\R}f(t)g(t)w^2(t)dt. \leqno{(3.1)} \]
[Proof] For $fw\in L^p(\R)$, $gw\in L^q(\R)$ such that $\ds \|gw\|_{L^q(\R)}=1, \int_\R g(t)w^2(t)dt=0$ and any fixed $a\in\R$, we have 
\begin{align*}
\ab{\int_\R f(t)g(t)w^2(t)dt} &= \ab{\int_\R(f(t)-a)g(t)w^2(t)dt} \tag{3.2} \\
&\le \norm{(f-a)w}{L^p(\R)}\norm{gw}{L^q(\R)} \\
&= \norm{(f-a)w}{L^p(\R)}. 
\end{align*}
That is
\[ E_{p,0}(w;f) \ge \sup_{\narabe{ \|gw\|_{L^q(\R)}\le 1 \\ \int_{\R}g(t)w^2(t)dt=0 }}\int_{\R}f(t)g(t)w^2(t)dt. \leqno{(3.3)} \]
Next, we will proof  the opposite direction. We first show the case of $1< p <\infty$. For $fw\in L^p(\R)$ and $c\in \R$, then $(f-c)w\in L^p(\R)$. We put
\begin{align*} 
I(c) :&= \int_\R\sign|f(t)-c|^{\frac{p}{q}}w^{\frac{p}{q}+1}(t)dt. \\
&= \int_\R\frac{(f(t)-c)}{|f(t)-c|}|f(t)-c|^{\frac{p}{q}}w^{\frac{p}{q}+1}(t)dt.
\end{align*}
Because $I(c)$ is continuous for every $c\in\R$ and $\ds \lim_{c\to\infty}I(c)=-\infty, \lim_{c\to-\infty}I(c)=\infty$, so we can find a constant $c_0$ such that $I(c_0)=0$. So we put
\[ g(t):=\frac{(f(t)-c_0)|f(t)-c_0|^{\frac{p}{q}-1}w^{\frac{p}{q}-1}(t)}{\|(f-c_0)w\|_{L^p(\R)}^\frac{p}{q}}, \]
then
\begin{align*} 
\|gw\|_{L^q(\R)} &= \kko{\int_{\R}\frac{|f(t)-c_0|^pw^p(t)}{\|(f-c_0)w\|_{L^p(\R)}^p}dt}^{\frac{1}{p}}= 1, \\
\end{align*}
\begin{align*} 
\int_\R g(t)w^2(t)dt &= \int_\R\sign\frac{|f(t)-c_0|^{\frac{p}{q}}}{\|(f-c_0)w\|^{\frac{p}{q}}}w^{\frac{p}{q}+1}(t)dt \\
&= \frac{I(c_0)}{\|(f-c_0)w\|^{\frac{p}{q}}}=0
\end{align*}
and 
\begin{align*}
\int_\R f(t)g(t)w^2(t)dt &= \int_\R (f(t)-c_0)g(t)w^2(t)dt \\
&= \int_\R (f(t)-c_0)\frac{(f(t)-c_0)|f(t)-c_0|^{\frac{p}{q}-1}w^{\frac{p}{q}+1}}{\|(f-c_0)w\|_{L^p(\R)}^\frac{p}{q}}dt \\
&= \frac{1}{\|(f-c_0)w\|_{L^p(\R)}^\frac{p}{q}} \int_\R|f(t)-c_0|^{\frac{p}{q}+1}w^{\frac{p}{q}+1}dt\\
&= \frac{1}{\|(f-c_0)w\|_{L^p(\R)}^\frac{p}{q}} \|(f-c_0)w\|_{L^p(\R)}^p \\
&= \|(f-c_0)w\|_{L^p(\R)}.
\end{align*}
That is,
\[ E_{p,0}(w;f) \le \sup_{\narabe{ \|gw\|_{L^q(\R)}\le 1 \\ \int_{\R}g(t)w^2(t)dt=0 }}\int_{\R}f(t)g(t)w^2(t)dt. \leqno{(3.4)} \]
If $p=1$,  we put
\begin{align*} 
I(c) :&= \int_\R\frac{(f(t)-c)}{|f(t)-c|}w(t)dt \\
&= \int_\R\sign w(t)dt.
\end{align*}
and let
\[ c_1:=\sup_c\{I(c)>0\},\ \ \ \  c_2:=\inf_c\{I(c)<0\}.  \]
We note that $c_1\le c_2$. If $c_1<c_2$, we put $c_0:=(c_1+c_2)/2$, then $I(c_0)\le 0$ and $I(c_0)\ge 0$, that is $I(c_0)=0$. On the other hand, if $c_1=c_2=c_0$, we put
\[ \lim_{c\to c_0}I(c)=A,\ \ \ \ \lim_{c\to c_0}I(c)=B, \]
\begin{align*}
E_+ &:=\set{x\in\R}{f(x)-c_0>0}, \\
E_- &:=\set{x\in\R}{f(x)-c_0<0}, \\
E_0 &:=\set{x\in\R}{f(x)=c_0}.
\end{align*}
We may suppose $B<A$, then 
\[ \int_{E_0}w(t)dt>0. \]
We note that 
\[ A= \int_{E_+\cup E_0}w(t)dt-\int_{E_-}w(t)dt \]
and
\[ B= \int_{E_+}w(t)dt-\int_{E_0\cup E_-}w(t)dt. \]
So 
\[ A-B= 2\int_{E_0}w(t)dt. \]
Now we consider the case of $A\le \int_{E_0}w(t)dt$ (tha case of $-B\le \int_{E_0}w(t)dt$ is similarly). We put 
\begin{align*}
g(t) := \frac{1}{w(t)}\times\begin{cases} 1 & t\in E_+ \\ 1-A/(\int_{E_0}w(t)dt) & t\in E_0 \\ -1 & t\in E_- \end{cases}.
\end{align*}
Then 
\[ \|gw\|_{L^\infty(\R)}=1 \]
and 
\begin{align*}
\int_\R g(t)w^2(t)dt &= \int_{E_+}w(t)dt+\int_{E_-}(-1)w(t)dt+\int_{E_0}w(t)g(t)dt \\
&= A-\int_{E_0}w(t)dt +\int_{E_0}w(t)g(t)dt \\
&= A-\int_{E_0}\frac{A}{\int_{E_0}w(t)dt}w(t)dt=0,
\end{align*}
furthermore,
\begin{align*}
\int_\R f(t)g(t)w^2(t)dt &= \int_\R (f(t)-c_0)g(t)w^2(t)dt \\
&= \int_{E_+} (f(t)-c_0)w(t)dt-\int_{E_-} (f(t)-c_0)w(t)dt\\
&\ \ +\int_{E_0} (f(t)-c_0)g(t)w^2(t)dt\\
&= \int_{E_+} |f(t)-c_0)|w(t)dt+\int_{E_-} |f(t)-c_0)|w(t)dt\\
&\ \ +\int_{E_0} |f(t)-c_0)|w(t)dt\\
&= \|(f-c_0)w\|_{L^1(\R)},
\end{align*}
That is,
\[ E_{1,0}(w;f) \le \sup_{\narabe{ \|gw\|_{L^\infty(\R)}\le 1 \\ \int_{\R}g(t)w^2(t)dt=0 }}\int_{\R}f(t)g(t)w^2(t)dt. \leqno{(3.5)} \]
Next, we show the case of $p=\infty$. For $fw\in L^\infty(\R)$, we put
\[ I(c):=\einf (f-c_0)w +\esup (f-c_0)w. \]
As the previous case, there exists a constant $c_0$ such that $I(c_0)=0$. 
For any small $\ve>0$, let
\begin{align*}
F_+ &:=\set{x\in\R}{(f(x)-c_0)w(x)>\|(f-c_0)w\|_{L^\infty(\R)}-\ve}, \\
F_- &:=\set{x\in\R}{(f(x)-c_0)w(x)<-\|(f-c_0)w\|_{L^\infty(\R)}+\ve}.
\end{align*}
Then $F_+$ and $F_-$ are disjoint. Let 
\[  A:=\int_\R\chi_{F_+}(t)w^2(t)dt>0,\ \ B:=\int_\R\chi_{F_-}(t)w^2(t)dt>0\]
 and
\[ C:=\frac{1}{2A}\int_\R\chi_{F_+}(t)w(t)dt + \frac{1}{2B}\int_\R\chi_{F_-}(t)w(t)dt.  \]
Moreover,
\[ g(t):=\frac{1}{C}\kko{\frac{1}{2A}\chi_{F+}(t)-\frac{1}{2B}\chi_{F_-}(t)}. \]
Then
\begin{align*}
\|gw\|_{L^1(\R)} &= \int_\R\frac{1}{C}\ab{\frac{1}{2A}\chi_{F+}(t)-\frac{1}{2B}\chi_{F_-}(t)}w(t)dt \\
&= \frac{1}{2C}\kko{\frac{1}{A}\int_\R\chi_{F_+}(t)w(t)dt+\frac{1}{B}\int_\R\chi_{F_-}(t)w(t)dt} \\
&= 1
\end{align*}
and
\begin{align*}
\int_\R g(t)w^2(t)dt &= C\kko{\frac{1}{2A}\int_\R\chi_{F_+}w^2(t)dt-\frac{1}{2B}\int_\R\chi_{F_-}w^2(t)dt}=0.
\end{align*}
Moreover, 
\begin{align*}
\int_\R f(t)g(t)w^2(t)dt &= \int_\R (f(t)-c_0)g(t)w^2(t)dt \\
&= \int_\R C\kko{\frac{(f(t)-c_0)}{2A}\chi_{F_+}-\frac{(f(t)-c_0)}{2B}\chi_{F_-}}w^2(t)dt \\
&= \int_{F_+}\frac{C(f(t)-c_0)}{2A}w^2(t)dt+\int_{F_-}\frac{-C(f(t)-c_0)}{2B}w^2(t)dt \\
&\ge \int_{F_+}\frac{C(\|(f-c_0)w\|_{L^\infty(\R)}-\ve)}{2A}w(t)dt\\
&\ \ \ \ \ \ +\int_{F_-}\frac{-C(\|(f-c_0)w\|_{L^\infty(\R)}-\ve)}{2B}w(t)dt \\
&= \|(f-c_0)w\|_{L^\infty(\R)}-\ve.
\end{align*}
So we have 
\[ E_{\infty,0}(w;f) \le \sup_{\narabe{ \|gw\|_{L^1(\R)}\le 1 \\ \int_{\R}g(t)w^2(t)dt=0 }}\int_{\R}f(t)g(t)w^2(t)dt. \leqno{(3.6)} \]
By (3.4), (3.5) and (3.6), the proof of the opposite direction is completed 
\vspace{2ex}
\qw{\bf Lemma 3.2} Let $P\in\sa{P}_n$. Then we have
\[ \int_\R(f(t)-v_n(f)(t))P(t)w^2(t)dt=0. \]
[Proof] Let $j\ge n+1$. For $P\in\sa{P}_n$, we have $\ds P(t)=\sum_{k=0}^na_kp_k(t)$. Therefore, we denote that $\ds c_k:=\int_\R f(x)p_k(x)w^2(x)dx$,
\begin{align*}
&\int_\R(f(t)-s_j(f)(t))P(t)w^2(t)dt \\
&= \int_\R\left(f(t)-p(t)\sum_{k=0}^{j-1}\int_\R f(x)p_k(x)w^2(x)dx\right)P(t)w^2(t)dt \\
&= \int_\R\left(f(t)-\sum_{k=0}^{j-1}c_kp_k(t)\right)\sum_{l=0}^na_lp_l(t)w^2(t)dt \\
&= \sum_{l=0}^na_lc_l-\sum_{k=0}^{j-1}a_kc_k =0.
\end{align*}
Since $\ds v_n(f)=\frac{1}{n}\sum_{j=n+1}^{2n}s_j(f)$, we have the result.\ \ \ \ $\Box$\\\\
\vspace{2ex}
\qw{\bf Lemma 3.3.} For $x\in\R$ and $n\in\N$, we have
\[ E_{1,2n-1}(w;\chi_x)\le C\frac{a_n}{n}w(x). \leqno{(3.7)} \]
[Proof] Let $n\in\N$ be a fixed and integer $0\le k\le n+1$ be found so that $x\in (x_{k+1,n},x_{k,n}]$. Then there exist $P:=P_k\in\sa{P}_{2n-2}$ and $R:=R_k\in\sa{P}_{2n-2}$, such that
\[ R_k(t) \le \chi_{(-\infty,x]}(t) \le P_k(t),\ \ \ \ (t\in\R) \leqno{(3.8)} \]
and
\[ \int_\R(P_k(t)-R_k(t))w^2(t)dt\le \lambda_{k+1,n}+\lambda_{k,n}. \leqno{(3.9)} \]
(see [2, Corollary 1.2.6].). Let $|x|\le a_{n/3}$. First, we estimate $E_{1,n}(w^2;\chi_{(-\infty,x]})$.\ let $k$ be an positive integer such tha $x\in [x_{k+1,n},x_{k,n}]\subset [-a_n,a_n]$.\ Here, we note that $x_{1,n}<a_n$.\ There exist polynomials $P$ and $R$ satisfying (3.8) and (3.9), so that
\begin{align*}
&E_{1,2n-1}(w^2;\chi_{(-\infty,x]})\\
\le & \int_\R(P(t)-\chi_{(-\infty,x]}(t))w^2(t)dt + \int_\R(\chi_{(-\infty,x]}(t)-R(t))w^2(t)dt \\
\le & \lambda_{k+1,n}+\lambda_{k,n}. 
\end{align*}
By Lemma 2.3 (a), we have
\[ \lambda_{k,n}\sim\varphi_n(x_{k,n})w^2(x_{k,n}), \]
and by Lemma 2.2 (b), if $|x-x_{k,n}|\le C\varphi_n(x)$, then $w(x)\sim w(x_{k,n})$. We also have $\varphi_n(x_{k,n})\sim\varphi_n(x)$ by Lemma 2.2 (a). Therefore
\[ \lambda_{k,n}\sim C\varphi_n(x)w^2(x). \]
Similarly, from $\varphi_n(x_{k+1},n)\sim\varphi_n(x_{k,n})\sim\varphi_n(x)$ and $|x-x_{k+1,n}|\le C\varphi_n(x)$, we have
\[ \lambda_{k+1,n}\sim C\varphi_n(x)w^2(x). \]
Conseqently
\[ E_{1,2n-1}(w^2;\chi_{(-\infty,x]})\le C\varphi_n(x)w^2(x)\sim\frac{a_n}{n}\sqrt{1-\frac{a_n}{n}}w^2(x)\le C\frac{a_n}{n}w^2(x). \]
Now, we estimate $E_{1,2n-1}(w;\chi_{x})$. Since the MRS number for the weight $w^2$ equals $a_{2n}$, 
\begin{align*}
E_{1,2n-1}(w;\chi_{x})=E_{1,2n-1}((w^{1/2})^2;\chi_{(-\infty,x]})\le C\frac{a_{2n}}{n}w(x)\le C\frac{a_n}{n}w(x).
\end{align*}
Next, we consider $|x|\ge a_{n/3}$.\ By Lemma 2.1 (2), we have
\[ Q^\prime(x)\ge Q^\prime(a_{n/3})\ge C\frac{(n/3)\sqrt{T(a_{n/3})}}{a_{n/3}}\ge C\frac{n\sqrt{T(a_{n})}}{a_{n}}. \]
Suppose $x>0$,
\begin{align*}
E_{1,2n-1}(w;\chi_{(-\infty,x]}) &\le E_{1,0}(w;\chi_{(-\infty,x]}) \\
&\le \int_\R(1-\chi_{(-\infty,x]}(t))w(t)dt = \int_x^\infty \exp(-Q(t))dt \\
&\le \frac{-1}{Q^\prime(x)}\int_x^\infty(-Q^\prime(t))\exp(-Q(t))dt = \frac{w(x)}{Q^\prime(x)} \\
&\le C\frac{a_n}{n\sqrt{T(a_n)}}w(x)\le C\frac{a_n}{n}w(x).\ \ \ \ \Box
\end{align*}
\vspace{2ex}\\
For $hw\in L^p(\R)$, we set
\[ I(h)(t):=\frac{1}{w^2(t)}\int_t^\infty h(u)w^2(u)du. \]
In following Lemma, weight $w$ is Erd{\H{o}}s type. 
\vspace{2ex}
\qw{\bf Lemma 3.4.} (1) If $hw\in L^p(\R)$ and
\[ \int_\R h(t)w^2(t)dt=0,  \leqno{(3.10)}\]
then
\[ \|I^\prime(h)w\|_{L^p(\R)}\le C\|hw\|_{L^p(\R)} \leqno{(3.11)} \]
Moreover, if $g$ is absolutely continuous in each compact interval and $wg^\prime \in L^q(\R)$, then 
\[ \int_\R g(x)h(x)w^2(x)dx = \int_\R g^\prime(t)I(h)(t)w^2(t)dt. \leqno{(3.12)} \]
(2) If $hw\in L^\infty(\R)$, and
\[ \int_\R h(t)P(t)w^2(t)dt=0,\ \ \ \ P\in\sa{P}_n, \leqno{(3.13)} \]
then
\[ \|I(h)w\|_{L^\infty(\R)}\le C\frac{a_n}{n}\|hw\|_{L^\infty(\R)}. \leqno{(3.14)} \]
[Proof] We see that
\begin{align*}
I^\prime(h)(t)w(t) &= \frac{2Q^\prime(t)}{w(t)}\int_t^\infty w^2(u)h(u)du-h(t)w(t) \\
&= 2Q^\prime(t)w(t)I(h)(t)-h(t)w(t).
\end{align*}
So, in order to prove (3.11), it is enough to show that
\[ \|Q^\prime I(h)w\|_{L^p(\R)}\le C\|hw\|_{L^p(\R)}. \leqno{(3.15)} \]
First we estimate for $p=\infty$. For $t>0$, we have
\[ \int_t^\infty w(x)dx \le \int_t^\infty w(x)\kko{1+\frac{Q^{\prime\prime}(x)}{{Q^\prime}^2(x)}}dx. \]
Since,
\[ \frac{d}{dx}\kko{\frac{w(x)}{Q^\prime(x)}}=-w(x)\kko{1+\frac{Q^{\prime\prime}(x)}{{Q^\prime}^2(x)}}, \]
we have
\[ \frac{w(t)}{Q^\prime(t)}=\int_t^\infty w(x)\kko{1+\frac{Q^{\prime\prime}(x)}{{Q^\prime}^2(x)}dx}\ge \int_t^\infty w(x)dx. \]
So, for all $t>0$, we have
\[ \frac{Q^\prime(t)}{w(t)}\int_t^\infty w(x)dx \le 1. \leqno{(3.16)} \]
Therefore
\begin{align*}
|Q^\prime(t)I(h)(t)w(t)| &= \ab{\frac{Q^\prime(t)}{w(t)}\int_t^\infty h(x)w^2(x)dx} \\
&\le \|hw\|_{L^\infty(\R)}\ab{\frac{Q^\prime(t)}{w(t)}\int_t^\infty w(x)dx} \\
&\le \|hw\|_{L^\infty(\R)}. \tag{3.17}
\end{align*}
If $t<0$, (3.16) gives
\begin{align*}
\ab{\frac{Q^{\prime}(t)}{w(t)}\int_{-\infty}^t w(x)dx} &= \frac{Q^\prime(|t|)}{w(|t|)}\int_{|t|}^\infty w(x)dx \\
&\le 1. \tag{3.18}
\end{align*}
We may rewrite the condition (3.10) in the form
\[ \int_t^\infty h(x)w^2(x)dx = -\int_{-\infty}^t h(x)w^2(x)dx, \]
so by (3.18), for $t<0$,
\begin{align*}
\ab{\frac{Q^\prime(t)}{w(t)}\int_t^\infty h(x)w^2(x)dx} &= \ab{\frac{Q^\prime(t)}{w(t)}\int_{-\infty}^t h(x)w^2(x)dx} \\
&\le \|hw\|_{L^\infty(\R)}\ab{\frac{Q^\prime(t)}{w(t)}\int_{|t|}^\infty w(x)dx} \\
&\le \|hw\|_{L^\infty(\R)}.
\end{align*}
Therefore, hor all $h$ satisfying (3.10), we have shown that
\[ \|Q^\prime I(h)w\|_{L^\infty(\R)}\le \|hw\|_{L^\infty(\R)}. \leqno{(3.19)} \]
Next, we estimate for $p=1$. Let $H(x):=|h(x)|+|h(-x)|$, we have
\begin{align*}
 &\|Q^\prime I(h)w\|_{L^1(\R)} = \int_\R\ab{\frac{Q^\prime(t)}{w(t)}\int_t^\infty h(x)w^2(x)dx}dt \\
&= \int_{0}^{\infty}\ab{\frac{Q^\prime(t)}{w(t)}\int_t^\infty h(x)w^2(x)dx}dt + \int_{-\infty}^{0}\ab{\frac{Q^\prime(t)}{w(t)}\int_t^\infty h(x)w^2(x)dx}dt \\
&= \int_{0}^{\infty}\ab{\frac{Q^\prime(t)}{w(t)}\int_t^\infty h(x)w^2(x)dx}dt + \int_{-\infty}^{0}\ab{\frac{Q^\prime(t)}{w(t)}\int_{-\infty}^t h(x)w^2(x)dx}dt \\
&= \int_{0}^{\infty}\ab{\frac{Q^\prime(t)}{w(t)}\int_t^\infty h(x)w^2(x)dx}dt + \int_{0}^{\infty}\ab{\frac{Q^\prime(s)}{w(s)}\int_{-\infty}^{-s} h(x)w^2(x)dx}ds \\
&\le \int_0^\infty \frac{Q^\prime(t)}{w(t)}\kko{\int_t^\infty H(x)w^2(x)dx}dt = \int_0^\infty H(x)w^2(x)\kko{\int_0^x\frac{Q^\prime(t)}{w(t)}dt}dx \\
&= \int_0^\infty H(x)w^2(x)\kko{\frac{1}{w(x)}-1}dx \le \int_0^\infty H(x)w^2(x)\frac{1}{w(x)}dx \\
&= \int_\R |h(x)|w(x)dx = \|hw\|_{L^1(\R)}.\\ \tag{3.20}
\end{align*}
Now, let $fw\in L^p(\R) (1\le p\le\infty)$. We define a linear operator 
\[ F(f):=\frac{\ds\int_\R f(x)w^2(x)dx}{\ds\int_\R w^2(x)dx}. \]
Then we have 
\[ |F(f)|\le C\|fw\|_{L^p(\R)} \tag{3.21} \]
for $p=1, \infty$. And so
\[ \|(f-F(f))w\|_{L^p(\R)}\le C\|fw\|_{L^p(\R)} \tag{3.22} \]
for $p=1, \infty$, that is $(f-F(f))w\in L^p(\R)$ for $p=1, \infty$. Next, we have
\begin{align*}
\int_\R(f(t)-F(f))w^2(t)dt &= \int_\R f(t)w^2(t)dt - F(f)\int_\R w^2(t)dt \\
&= \int_\R f(t)w^2(t)dt - \int_\R f(x)w^2(x)dx =0,
\end{align*}
that is $h:=f-F(f)$ satisfies (3.10). Let $p=1, \infty$ and $fw\in L^p(\R)$, we consider the operator
\[ U(f):=Q^\prime I(f-F(f)). \]
Then, from (3.17), (3.20) and (3.22), we conclude that
\[ \|U(f)w\|_{L^p(\R)}\le C\|(f-F(f))w\|_{L^p(\R)}\le C\|fw\|_{L^p(\R)} \tag{3.23} \]
By Riesz-Thorin interpolation theorem, the estimate (2.28) holds for all $1\le p\le\infty$. Therefore
\[ \|U(f)w\|_{L^p(\R)}=\|Q^\prime I(f-F(f))w\|_{L^p(\R)}\le C\|fw\|_{L^p(\R)}. \tag{3.24} \]
Now we see that if $h$ satisfies (3.10), then we have $F(h)=0$.\ So applying (3.24) to $f=h$, we conclude
\[ \|Q^\prime I(h)w\|_{L^p(\R)}\le C\|hw\|_{L^p(\R)} \]
for all $1\le p\le\infty$.\ This completes the proof of (3.15), that is (3.11).

Next, for $1<p<\infty$, in the same way as (3.18), for $t<0$, we have
\[ |qQ^\prime(t)|\int_{-\infty}^t w^q(x)dx\le w^q(t) \]
and since $w$ is Erd{\H{o}}s type, we have 
\[ \frac{1}{|Q^\prime(t)|}\le C\frac{1}{|t|^{2q/p}} \]
holds (see [4, (2.4)]). Now let $g^\prime w\in L^q(\R)$, then $gw\in L^q(\R)$ (cf. [4, Theorem 6]), hence
\begin{align*}
&\int_{-\infty}^{-1}|g^\prime(t)|\kko{\int_{-\infty}^t|h(x)|w^2(x)dx}dt \\
&\le \int_{-\infty}^{-1}|g^\prime(t)|\kko{\int_{-\infty}^t(|h(x)|w(x))^pdx}^{1/p}\kko{\int_{-\infty}^tw^q(x)dx}^{1/q}dt \\
&\le \|hw\|_{L^p(\R)}\int_{-\infty}^{-1}|g^\prime(t)|\kko{\frac{w^q(t)}{Q^\prime(t)}}^{1/q}dt\\
&= \|hw\|_{L^p(\R)}\int_{-\infty}^{-1}|g^\prime(t)|w(t)\frac{1}{(Q^\prime(t))^{1/q}}dt\\
&\le \|hw\|_{L^p(\R)}\|g^\prime w\|_{L^q(\R)}\kko{\int_{-\infty}^{-1}\frac{1}{t^2}dt}^{1/p} <\infty.
\end{align*}
If $p=\infty$, then 
\begin{align*}
&\int_{-\infty}^{-1}|g^\prime(t)|\kko{\int_{-\infty}^t|h(x)|w^2(x)dx}dt \\
&\le \int_{-\infty}^{-1}\kko{|g^\prime(t)|\|hw\|_{L^\infty(\R)}\int_{-\infty}^{t}w(x) dx}dt\\
&\le \|hw\|_{L^\infty(\R)}\int_{-\infty}^{-1}\frac{w(t)}{|Q^\prime(t)|}|g^\prime(t)|dt \\
&\le C\|hw\|_{L^\infty(\R)}\int_{-\infty}^{-1}|g^\prime(t)|w(t)dt<\infty.
\end{align*}
Finally, we consider $p=1$, then 
\begin{align*}
&\int_{-\infty}^{-1}|g^\prime(t)|\kko{\int_{-\infty}^t|h(x)|w^2(x)dx}dt \\
&=\int_{-\infty}^{-1}|h(x)|w^2(x)\kko{\int_{x}^{-1}|g^\prime(t)|dt}dx \\
&=\int_{-\infty}^{-1}|h(x)|w(x)\kko{\int_{x}^{-1}|g^\prime(t)|w(t)\frac{w(x)}{w(t)}dt}dx. \\
\end{align*}
We note that
\[ \int_x^{-1}\frac{w(x)}{w(t)}dt=w(x)\int_x^{-1}\frac{Q^\prime(t)}{w(t)}\frac{1}{Q^\prime(t)}dt\le \frac{w(x)}{C}\int_x^{-1}\frac{Q^\prime(x)}{w(t)}dt<\frac{1}{C}. \]
holds. So 
\[ \int_{-\infty}^{-1}|g^\prime(t)|\kko{\int_{-\infty}^t|h(x)|w^2(x)dx}dt<\infty. \]
And
\[ \int_{-1}^0|g^\prime(t)|\kko{\int_{-\infty}^t|h(x)|w^2(x)dx}dt<\infty. \]
By (3.10), we may assume thet $g(0)=0$ in (3.12).\ Then
\begin{align*}
&\int_{-\infty}^\infty g(x)h(x)w^2(x)dx \\
&= -\int_{-\infty}^0h(x)w^2(x)\kko{\int_x^0g^\prime(t)dt}dx \\
&= -\int_{-\infty}^0\kko{\int_{-\infty}^0g^\prime(t)\chi_{[x,0]}(t)dt}h(x)w^2(x)dx \\
&= -\int_{-\infty}^0\kko{\int_{-\infty}^0h(x)\chi_{(-\infty,t]}(x)w^2(x)dx}g^\prime(t)dt \\
&= -\int_{-\infty}^0g^\prime(t)\kko{\int_{-\infty}^th(x)w^2(x)dx}dt \\
&= \int_{-\infty}^0g^\prime(t)\kko{\int_t^\infty h(x)w^2(x)dx}dt \\
&= \int_{-\infty}^0g^\prime(t)I(h)(t)w^2(t)dt.
\end{align*}
Similarly,
\[ \int_0^\infty g(x)h(x)w^2(x)dx=\int_0^\infty g^\prime(t)I(h,t)w^2(t)dt. \]
Hence, (3.12) is proved.

Moreover, let $h$ satisfies (3.13), we get for $t\in\R$ and $P_{\chi,t}\in\sa{P}_{2n-1}$ is the best approximation for $\chi_{[-\infty,t]}$ in $L^1(\R)$. Then 
\begin{align*}
|I(h)(t)| &= \ab{\frac{1}{w^2(t)}\int_t^\infty h(x)w^2(x)dx} \\
&= \ab{\frac{1}{w^2(t)}\int_{-\infty}^t h(x)w^2(x)dx} \\
&= \ab{\frac{1}{w^2(t)}\int_{\R}\chi_{(-\infty,t]}(x)h(x)w^2(x)dx} \\
&= \ab{\frac{1}{w^2(t)}\int_{\R}(\chi_{(-\infty,t]}(x)-P_{\chi,t}(x))h(x)w^2(x)dx} \\
&\le \frac{\|hw\|_{L^\infty(\R)}}{w^2(t)}\int_\R|\chi_{(-\infty,t]}(x)-P_{\chi,t}(x)|w(x)dx \\
&\le \frac{\|hw\|_{L^\infty(\R)}}{w^2(t)} E_{1,2n-1}(w;\chi_{(-\infty,t]})\\
&\le C\frac{1}{w(t)}\frac{a_n}{n}\|hw\|_{L^\infty(\R)},
\end{align*}
By Lemma 3.3. Therefore
\[ \|I(h,t)w\|_{L^\infty(\R)}\le C\frac{a_n}{n}\|hw\|_{L^\infty(\R)}.\ \ \ \ \Box \]
\vspace{1ex}
\qw{\bf Lemma 3.5.} [4, Theorem 1]. Let $g$ be absolutely continuous in each compact interval and $g^\prime w\in L^p(\R)$. Then we have
\[ E_{p,n}(w;g)\le C\frac{a_n}{n}\|g^\prime w\|_{L^p(\R)}, \]
and equivalently,
\[ E_{p,n}(w;g)\le C\frac{a_n}{n}E_{p,n-1}(w;g^\prime). \]
\vspace{1ex}
\qw{\bf Lemma 3.6.} [5. Theorem 1.1]. Let $1\le p \le \infty$ and $T(a_n)\le c(n/a_n)^{2/3}$ for some $c>0$. If $T^{1/4}fw\in L^p(\R)$, then
\[ \|(f-v_n(f))w\|_{L^p(\R)} \le CT^{1/4}(a_n)E_{p,n}(w;f). \]
\vspace{2ex}
\qw{\bf Lemma 3.7.} Let $g$ be absolutely continuous and $g^\prime\in L^{p}(\R)$, $T(a_n)\le c(n/a_n)^{2/3}$ for some $c>0$.. Then for $n\in\N$ there exists a polynomial $V_n(g)\in\sa{P}_{2n}$ such that $V_n^\prime=v_n(g^\prime)$ and
\[ \|(g-V_n)w\|_{L^p(\R)}\le C\frac{a_n}{n}T^{1/4}(a_n)E_{p,n}(w;g^\prime). \leqno{(3.25)} \]
[Proof] Without loss of generality, we may assume the $g(0)=0$. Let
\[ G(x):=\int_0^x(g^\prime(t)-v_n(g^\prime)(t))dt, \]
and take a constant $a$ such that
\[ \|(G-a)w\|_{L^p(\R)}\le 2E_{p,0}(w;G). \] 
Then
\begin{align*}
V_n(x) :&= a+\int_0^x v_n(g^\prime)(t)dt \\
&= a+g(x)-G(x) \in \sa{P}_{2n}, \tag{3.26}
\end{align*}
and
\[ \|(g-V_n)w\|_{L^p(\R)}=\|(G-a)w\|_{L^p(\R)}\le 2E_{p,0}(w;G). \leqno{(3.27)} \]
By the duality principle and Lemma 3.1, we find $h$ such thae $hw\in L^q(\R)$ and satisfies (3.1), $\|hw\|_{L^q(\R)}=1$ and by Lamma 3.2, in the same way as changing the order of integration in proof of (3.12), 
\begin{align*}
& \|(g-V_n(g))w\|_{L^p(\R)} \le 2E_{p,0}(w;G) \\
&\le 4\ab{\int_\R G(x)h(x)w^2(x)dx} \\
&= 4\int_\R \kko{\int_0^xg^\prime(t)-v_n(g^\prime)(t)dt}h(x)w^2(x)dx \\
&= 4\int_\R(g^\prime(t)-v_n(g^\prime)(t))\kko{\int_t^\infty h(x)w^2(x)dx}dt \\
&= 4\int_\R(g^\prime(t)-v_n(g^\prime)(t))I(h)(t)w^2(t)dt \\
&= 4\int_\R(g^\prime(t)-v_n(g^\prime)(t))(I(h)(t)-P_{I(h)})w^2(t)dt \\
\end{align*}
holds(cf. [4. Theorem 6]).  Where $P_{I(h)}\in\sa{P}_n$ is the best approximation polynomial for $I(h)$ with $L^q(\R)$.\ From Lemma 3.5 with $q$ and Lemma 3.4 (3.11) with $q$, we have
\begin{align*}
E_{q,n}(w;I(h)) &\le C\frac{a_n}{n}\|I^\prime(h) w\|_{L^q(\R)} \le C\frac{a_n}{n}\|hw\|_{L^q(\R)} = C\frac{a_n}{n}. \tag{3.28}
\end{align*}
Since by Lemma 3.6,
\begin{align*}
\|(g^\prime-v_n(g^\prime)w\|_{L^p(\R)} &\le CT^{1/4}(a_n)E_{p,n}(w;g^\prime). \tag{3.29}
\end{align*}
From (3.28) and (3.29), we have the result.\ \ \ \ $\Box$.\\

\section{{Proof of Theorem 1.1}}
\qw{\bf Theorem 4.1.} Let $w\in\sa{F}(C^2+)$ is Erd{\H{o}}s type and suppose that $T(a_n)\le c({n}/{a_n})^{2/3}$ holds some constant $c>0$. Then there exists a constant $C\ge 1$ such that for any absolutely continuous function $f$ with $f^\prime w \in L^q(\R)$, if
\[ \|(f-P)w\|_{L^p(\R)}\le \eta \]
holds for some $n\in\N$, $P\in\sa{P}_n$ and $\eta>0$, then
\[ \|w(f^\prime-P^\prime)\|_{L^p(\R)}\le C\kko{T^{3/4}(a_n)E_{p,n}(w;f^\prime)+\frac{n}{a_n}T^{1/2}(a_n)\eta}. \leqno{(4.1)} \]
$[$Proof$]$   By [2, Theorem 1.15 and Corollary 1.16], we see  
\[ \|P^\prime w\|_{L^p(\R)} \le C\frac{nT^{1/2}(a_n)}{a_n}\|Pw\|_{L^p(\R)} \]
  and by Proposition 3.3, we have  
\[ \|(f-V_n)w\|_{L^p(\R)} \le C\frac{a_{n}}{n}T^{1/4}(a_{n})E_{p,n}(w;f^\prime) \]
where $V_n^\prime=v_n(f^\prime)$.   Also by (3.29), we have  
\begin{align*}
\|(f^\prime-v_n(f^\prime))w\|_{L^p(\R)} &\le CT^{1/4}(a_n)E_{p,n}(w;f^\prime). \tag{4.2}
\end{align*}
Hence  
\begin{align*}
&\|(f^\prime-P^\prime)w\|_{L^p(\R)} \le \|(f^\prime-v_n(f^\prime))w\|_{L^p(\R)}+\|(V_n^\prime-P^\prime)w\|_{L^p(\R)} \\
&\le  CT^{1/4}(a_n)E_{p,n}(w;f^\prime)+C\frac{2n}{a_{2n}}T^{1/2}(a_{2n})\|(V_n-P)w\|_{L^p(\R)} \\
&\le  CT^{1/4}(a_n)E_{p,n}(w;f^\prime) \\
&\ \ \ \ \ \ +C\frac{n}{a_{n}}T^{1/2}(a_{n})(\|(f-V_n)w\|_{L^p(\R)}+\|(f-P)w\|_{L^p(\R)})\\
&\le C\kko{T^{1/4}(a_n)E_{p,n}(w;f^\prime)+T^{3/4}(a_n)E_{p,n}(w;f^\prime)+\frac{n}{a_n}T^{1/2}(a_n)\eta},
\end{align*}
  which shows (4.1).\ \ \ \ $\Box$\\
\qw{[Proof of Theorem 1.1]} In (4.1), we take $\eta=C_1T^{1/4}(a_n)E_{p,n}(w;f)$.  
Then since Lemma 3.5, we have  
\begin{align*}
&\|(f^\prime-P^\prime)w\|_{L^p(\R)} \\
&\le C_2\kko{T^{3/4}(a_n)E_{p,n}(w;f^\prime)+C_1\frac{n}{a_n}T^{3/4}(a_n)E_{p,n}(w;f)}\\
&\le C_2\kko{T^{3/4}(a_n)E_{p,n-1}(w;f^\prime)+C_1T^{3/4}(a_n)E_{p,n-1}(w;f^\prime)} \\
&=C_2(1+C_1)T^{3/4}(a_n)E_{p,n-1}(w;f^\prime).
\end{align*}
  The desired result follows from above.\ \ \ \ $\Box$\\

Using Theorem 4.1 inductively, we have  an estimate in the case of higher order derivative.

\qw{\bf Corollary 4.2.} Let $w=\exp(-Q)$ be the same as in Proposition 4.1. and $j\in\N$.  For any $C_1\ge 1$, $1\le p\le \infty$ and absolutely continuous function $f^{j-1}$ with $f^(j) w\in L^p(\R)$, if
\[ \|(f-P)w\|_{L^p(\R)}\le C_1E_{p,n}(w;f) , \]
holds for some $n>j$ and $P\in\sa{P}_n$, then there exists constant $C_2\ge 1$
\[ \|w(f^{(i)}-P^{(i)})\|_{L^p(\R)}\le C_2(1+C_1)T^{(2i+1)/4}(a_n)E_{p,n-i}(w;f^{(i)}) \]
holds for all $i=1, 2, \cdots , j$.\\
$[$Proof $]$   We note if $f^\prime w\in L^p(\R)$ implies $fw\in L^p(\R)$,   it follows from $f^{(j)} w\in L^p(\R)$ that $f^{(i)} w\in L^p(\R)$ for all $0\le i \le j-1$.   We see 
\[ \|(f^\prime-P^\prime)w\|_{L^p(\R)} \le C_2(1+C_1)T^{3/4}(a_n)E_{p,n-1}(w;f^\prime). \]
  Applying Theorem 4.1 to $n-1$, $f^\prime,\ P^\prime$ and 
\begin{align*}
\eta :&= C_2(1+C_1)T^{3/4}(a_n)E_{p,n-1}(w;f^\prime) \\
&\le C_2(1+C_1)T^{3/4}(a_n)\frac{a_{n-1}}{n-1}E_{p,n-2}(w;f^{\prime\prime}),
\end{align*}
we have  
\begin{align*}
&\|(f^{\prime\prime}-P^{\prime\prime})w\|_{L^p(\R)}\\
&\le C_2\kko{T^{3/4}(a_{n-1})E_{p,n-1}(w;f^{\prime\prime})+\frac{n-1}{a_{n-1}}T^{1/2}(a_{n-1})\eta}\\
&\le C_2(1+C_1)T^{5/4}(a_n)E_{p,n-2}(w;f^{\prime\prime}),
\end{align*}
  This shows the case $i=2$.   Repeating this process, we have the result for $1\le i \le j$.\ \ \ \ $\Box$ 

\mbox{}\\
{\footnotesize{
\qw Kentaro Itoh\\
Department of Mathematics\\
Meijo University\\
Tenpaku-ku, Nagoya 468-8502\\
Aichi Japan\\
{\texttt{133451501@ccalumni.meijo-u.ac.jp}}
\vspace{1ex}
\qw Ryozi Sakai\\
Department of Mathematics\\
Meijo University\\
Tenpaku-ku, Nagoya 468-8502\\
Aichi Japan\\
{\texttt{ryozi@crest.ocn.ne.jp}}
\vspace{1ex}
\qw Noriaki Suzuki\\
Department of Mathematics\\
Meijo University\\
Tenpaku-ku, Nagoya 468-8502\\
Aichi Japan\\
{\texttt{suzukin@meijo-u.ac.jp}}
}}
\end{document}